\documentclass[11pt]{amsart}

\usepackage{mathpple}
\usepackage{amsmath,amsfonts, amscd,graphicx}
\usepackage{amsrefs}
\usepackage{amsthm}
\usepackage[all]{xy}
\usepackage{color}

\usepackage{hyperref}

\usepackage{mathtools}

\usepackage{appendix}

\numberwithin{equation}{section}

\newcommand{\C}{\mathbb{C}}

\newcommand{\bracket}[1]{\left(#1\right)}

\newcommand{\mbf}{\mathbf}

\newcommand{\mc}{\mathcal}

\newcommand{\into}{\hookrightarrow}

\usepackage{amssymb}
\usepackage{verbatim}
\usepackage{mathrsfs}

\newcommand{\LD}{\langle}
\newcommand{\RD}{\rangle}

\DeclareMathOperator{\Jac}{Jac}

\theoremstyle{plain}
\newtheorem{thm}{Theorem}[section]

\newtheorem{thm-defn}{Theorem/Definition}[section]

\newtheorem{lem-defn}[thm]{Lemma/Definition}
\newtheorem{prop}[thm]{Proposition}

\newtheorem{prop-defn}[thm]{Proposition-Definition}

\newtheorem{thm-alg}[thm]{Theorem/Algorithm}

\allowdisplaybreaks[4]  

\begin{document}
\setlength{\unitlength}{1mm}

 \title{From primitive form to mirror symmetry}

 \author{Kyoji Saito}

  \address{Kavli Institute for the Physics and Mathematics of the Universe (WPI),
Todai Institutes for Advanced Study, The University of Tokyo,
5-1-5 Kashiwa-no-Ha,Kashiwa City, Chiba 277-8583, Japan}
\email{kyoji.saito@ipmu.jp}
\thanks{
 ${}^1$ Present note is worked out with the help of the coauthors Changzheng Li, Si Li and Yefeng Shen, to whom the author expresses his deep gratitudes.}


  \maketitle
\thispagestyle{plain}

\begin{abstract} This is a report on the recent joint work \cite{LLSS} on LG-LG mirror symmetry for the 14 exceptional unimodular singularities.${}^{1}$
\end{abstract}


\addtocontents{toc}{\protect\setcounter{tocdepth}{1}}

\section{Introduction}
In the early 80's, the author introduced the theory of primitive forms
\cite{Saito-lecutures, Saito-residue, Saito-universal,Saito-primitive},
which studies the period integrals of a primitive form over cycles which are vanishing 
to isolated critical points of a function.  One important consequence
of the theory  is that a primitive form induces a {\it  flat structure}   \cite{Saito-primitive} on the deformation parameter space of the function, including the {\it flat metric}, the {\it flat coordinate system}  and the {\it potential function} (which is later called the prepotential).
Later on, in the early 90's, Dubrovin \cite{D} studied the 2D topological field theory of genus zero curves and found the same structure, which is axiomatized to, so-called, the  \emph{Frobenius manifolds structures.} Lots of Frobenius manifolds structures are found. They include the examples constructed on the orbit spaces of Coxeter groups \cite{Saito-a},\cite{SYS},\cite{DZ},\cite{ShiraishiTakahashi} (which contain the first cases found before the primitive form theory), the one constructed by primitive forms \cite{Saito-primitive},\cite{D-S-I},\cite{D-S-II},\cite{IST},\cite{Takahashi}, the quantum cohomology rings,  Barannikov-Kontsevich construction \cite{BK} using polyvector fields of a  Calabi-Yau manifold, and the FJRW-thoery (the A-model of quantum singularity theory) \cite{FJR, FJR2}. 

The Gromov-Witten theory \cite{KM} counts pseudoholomorphic curves in a given symplectic manifold. Its application to  the symplectic structure on a K\"ahler manifold was extensively studied from a view point of the mirror symmetry. Here the mirror symmetry is one of the dualities in physics and had a strong impact on mathematics \cite{HV,mirror-book}. Namely, it asserts that certain data counted from symplectic geometry (the A-model side) should be equivalent to that from the complex structure of the \emph{mirror manifold} (the B-model side). There are several different levels of formulation of  the mirror symmetry such as  the categorical level \cite{Kontsevich,FOOO}, the geometric level \cite{SYZ}, or the equivalence of the genus zero theory on the A-model side with the variation of Hodge structures on the B-model side \cite{LLY, G3}. 

Primitive forms are about universal deformations $F$ of functions, giving flat structures on the deformation spaces. Hence, the theory is relevant in the complex geometric (B-model) aspects of N=(2,2) supersymmetric Landau-Ginzburg (LG) theory with the superpotential $F$. However,  this pattern of the mirror symmetry was not mathematically rigorously worked out until recently. This is because of  (1) lack of mathematical theory of A-model LG-theory at that early time, and (2) the difficulty of calculating primitive forms until recently, where explicit expressions of primitive forms were known only for weighted homogeneous polynomials of central charge less than or equal to 1. Both difficulties were resolved as follows. 

In around 2007, Fan, Jarvis and Ruan constructed a so-called \emph{quantum singularity  theory} by counting virtual cycles associated with a weighted homogeneous polynomial, whose potential (generating series) gives again a Frobenius manifold structure on the so-called \emph{FJRW state space} \cite{FJR, FJR2}. This is considered as an A-model Landau-Ginzburg theory. They immediately realized that such Frobenius manifold for an ADE-polynomial $W$ is actually isomorphic to the Frobenius manifold arising from the primitive form (i.e.\ B-side) of another ADE-polynomial $W^T$. The  superpotential polynomial $W^T$ on B-side is obtained by the transposition of exponents of monomials in the polynomial $W$ on A-side  \cite{BH,BHe,K}, which  is later called \emph{Berglund-H\"ubsch-Krawitz mirror}.
 As an application of this mirror theorem, they solved the so-called \emph{generalized Witten conjecture}, which says that the generating functions arising from the Landau-Ginzburg model for ADE-singularities should be governed by some ADE-integrable hierarchies.  Also a similar observation for simply elliptic singularities \cite{Saito-simplyElliptic} was achieved \cite{KS,MR,MS}. We remark that the relationship between FJRW theory and Gromov-Witten theory (both are in A-model side) is studied under the name of LG/CY-correspondence, for which one is referred to \cite{CIR,CR,CR2,KS,MR,R}.
 
On the other hand, in 2013, jointly with Changzheng Li and Si Li, the author came to a new perturbative construction of primitive forms \cite{LLSaito}, where Birkhoff decomposition theorem used in the original formulation \cite{Saito-primitive} was replaced by the asymptotic expansion of oscillatory integrals. This enables us to calculate primitive forms explicitly as a power series in an algorithmic way (at least for weighted homogeneous polynomials). 
With the perturbative approach, we can calculate further the flat coordinate system and the pre-potential function up to any finite order. This will be sufficient to determine the flat coordinate system and the pre-potential function with a help of WDVV-equations.

These two new developments thoroughly changed the view on the LG-LG mirror symmetry.  Namely, up to a choice of  primitive forms, one asks whether the pre-potential attached to FJRW theory for a weighted homogeneous polynomial could coincide with the prepotential associated to a primitive form for the mirror dual-polynomial.  Such a mirror symmetry  is called the Landau-Ginzburg to Landau-Ginzburg \emph{(LG-LG) mirror symmetry}. Its study has developped rapidly in the last years.  In the present note, we briefly introduce  the theories on both sides and the mirror map construction connecting them. Then, we confirm the LG-LG mirror symmetry for the 14 unimodular singularities, which are the first case of weighted homogeneous polynomials whose central charge exceeds 1.

\medskip
\noindent
{\bf Remark.}
We remark that the primitive form theory depends only on the
analytic equivalence class of the singularity of the function $W^T$,
although associated primitive forms may not be unique but form a family.
On the other hand, FJRW theory depends on the polynomial $W$ itself
together with a symmetry group of $W$. Hence, to achieve the mirror
symmetry  to FJRW theory on $W$, both the analytic equivalence class 
of $W^T$ and the choice of
the primitive form for $W^T$ depend on the choice of the polynomial $W$. We do not
yet understand this phenomenon conceptually (c.f.\  \cite{LLSS} {\bf Remark 4.9.} (2)).

\section{Primitive form theory.} The origin of a Landau-Ginzburg B-model (with respect to trivial group symmetry) at genus zero is the  theory of primitive forms   \cite{LLSaito, Saito-lecutures, Saito-residue, Saito-universal,Saito-primitive, Saito-Takahashi}.
The  starting data of the theory is a holomorphic function $
f:(X,\mathbf{0})\to (\C,{0})
$ defined on a Stein domain $X\subset \C^n$ with finite critical points. For our purpose on the LG-LG mirror symmetry, it is sufficient to consider
a     weighted homogeneous polynomial  $f=f(x_1, \cdots, x_n)$  with an isolated critical point at the origin $\mathbf{0}\in X=\C^n$,
$$
    f(\lambda^{q_1}x_1,\cdots, \lambda^{q_n}x_n)=\lambda  f(x_1,\cdots, x_n), \,\,\,\forall \lambda\in\C^*,
$$
  Here $\bracket{q_1,\cdots, q_n}$ in $\mathbb{Q}_{>0}^n$ are called the weights of the coordinates $(x_1, \cdots, x_n)$, and each weight $0<q_i\leq {1\over 2}$ is unique \cite{Saito-quasihomogeneous}.
In \cite{Saito-residue}, the author  introduced the  formal completion of the  Brieskorn lattice  together with a semi-infinite $z$-adic filtration by a formal variable $z$:
$$
 \hat{\mc{H}}_f^{(0)}:=\Omega^n_{X,\mbf{0}}[[z]]/(df+zd)\Omega^{n-1}_{X,\mbf{0}}[[z]],
$$
and constructed a  \emph{higher residue pairing}
$$
   K_f: \hat{\mc{H}}_f^{(0)}\otimes \hat{\mc{H}}_f^{(0)}\to z^{n}\C[[z]]
$$
which satisfies a number of properties, and plays a key role in the theory of primitive forms.
A universal unfolding of $f$ is given by
$$
  F:  (X\times S, \mathbf{0}\times \mathbf{0})\to  (\C, 0),\quad
 F(\mathbf{x}, \mathbf{s})=f(\mathbf{x})+\sum_{\alpha=1}^\mu s_\alpha \phi_\alpha, $$
where $\{\phi_1, \cdots, \phi_\mu\}\subset \C[\mathbf{x}]$ are weighted homogeneous polynomials representing an additive basis of the Jacobian algebra $\Jac(f)$, and $\mathbf{s}=\{s_\alpha\}_{\alpha=1,\cdots,\mu}$ parametrizes the deformation space $S\subset \C^\mu$. Using, so  called, Kodaira-Spencer map: $\sum_i a_i\partial_{s_i} \mapsto \sum_i a_i\phi_i$, the tangent bundle of $S$ is identified with the Jacobi ring of $F$, which gives a ring structure (Frobenius algebra structure) and a natural inner product $J$ (the first residue pairing) on the tangent bundle of $S$. 
There is a family version $\hat{\mc{H}}_F^{(0)}$ (resp. $K_F$) of $\hat{\mc{H}}_f^{(0)}$ (resp. $K_f$) with respect to the universal unfolding $F$.
 We remark that in the recent    work  \cite{LLSaito} by C. Li, S. Li and the author,  an alternate complex differential geometric approach to the module $\hat{\mc{H}}_F^{(0)}$
  is developed. Therein we give a simple construction of the higher residue pairing by using integration of compactly supported polyvector fields.

A primitive form is a section $\zeta\in \Gamma(S, \hat{\mc{H}}_F^{(0)})$,  represented by  a relative holomorphic volume form
$
\zeta=P(\mathbf{x},\mathbf{s})d^n\mathbf{x}$   ($d^n\mathbf{x}=dx_1\cdots dx_n$) 
on $X\times S$, satisfying the properties of \textit{primitivity, orthogonality, holonomicity, and homogeneity}, described by bilinear equations on $\zeta$ using the higher residue pairing $K_F$ together with  Gauss-Manin connection on $\hat{\mc{H}}_F^{(0)}$. 

Roughly speaking, the submodule of $ \hat{\mc{H}}_F^{(0)}$ consisting of the covariant differentiations of a primitive form by the tangent vectors of $S$ forms a splitting factor to the adic filtration on $ \hat{\mc{H}}_F^{(0)}$ defined by (the powers of) $z$, i.e. $\hat{\mc{H}}_F^{(0)} \simeq \mathcal{T}_S\ \oplus\  z\cdot \hat{\mc{H}}_F^{(0)}$. In this way, properties of the primitive form are transferred to the splitting factor i.e.\ to the tangent bundle of $S$, and, hence, the space $S$ obtains a  differential geometric structure, called the {\it flat structure} (= the {\it Frobenius manifold structure}) associated with $\zeta$. For instance, the orthogonality property of $\zeta$ gives a flat metric $J$ (i.e.\ the first residue pairing) on $S$. Then the flat section of the Levi-Civita  connection of that metric defines the flat coordinate system (see e.g. \cite{Saito-Takahashi} for more details).    

For  weighted homogeneous polynomials,  $\{\phi_\alpha d^n\mathbf{x}\}_\alpha \subset \hat{\mc{H}}_f^{(0)}$  is called a good basis if
the vector subspace $B=\mbox{Span}_\C\{\phi_\alpha d^n\mathbf{x}\}_\alpha$ satisfies $\mc K_f(B, B)\subset \C z^n$, where we note that the space $B$  is isomorphic to the Jacobi algebra $\Jac(f)$ as $\C$-vector spaces. 
One key step to construct a primitive form is that {\it the concept of the primitive forms is equivalent to the notion of   good section} \cite{Saito-primitive} (cf.\ \cite{Saito-existence}).  In order to show this, we need to extend  a good basis in  $ \hat{\mc{H}}_f^{(0)}$ to a "deformed good basis"  in the deformed module $ \hat{\mc{H}}_F^{(0)}$, where, in the proof, we use a classical analytic result known as {\it Birkhoff decomposition theorem}.
In  \cite{LLSaito}, we replaced the role of the Birkoff theorem by a multiplication of the "holomorphic part of the oscillatory integral factor" $e^{F-f\over z}:B\to B((z))[[\mathbf{s}]]$
 (here   the first copy of $B$ should be read off a subspace of the deformation $\hat{\mc{H}}_F^{(0)}$ of $\hat{\mc{H}}_f^{(0)}$), which is able to calculate in power series in the local coordinate perturbatively.
   Inspired from this, we obtain the following, which is a combination of several propositions  in section 3.2 of \cite{LLSS}.
   \begin{prop}\label{prop-pert}
 Given a good basis  $\{[\phi_\alpha d^n \mathbf{x}]\}_{\alpha=1}^\mu\subset \hat{\mc{H}}_f^{(0)}$,  there exists a unique  pair $(\zeta, \mc J)$ satisfying the following:
  $(1)\, \zeta \in B[[z]][[\mathbf{s}]], \quad (2)\, \mc J \in [d^n\mathbf{x}]+z^{-1}B[z^{-1}] [[\mathbf{s}]]\subset \hat{\mc{H}}_f[[s]],$ and
\begin{equation}\label{pert-eqn}
   e^{(F-f)/z} \zeta= \mc J. \tag{$\star$}
\end{equation}
 Moreover, we  embed  $z^{-1}\C[z^{-1}][[\mathbf{s}]]\into z^{-1}\C[[z^{-1}]][[\mathbf{s}]]$ and  decompose
$$
    \mc J=[d^nx]+\sum_{m=-1}^{-\infty} z^{m} \mc J_{m}, \quad \text{where}\ \mc J_m=\sum_\alpha \mc J_{m}^\alpha [\phi_\alpha d^n\mathbf{x}], \mc J_{m}^\alpha\in \C[[\mathbf{s}]].
$$
Then $\zeta$ gives a formal primitive form, and $\{\mc J_{-2}^\alpha\}$ give a formal Frobenius manifold structure on $S$ with flat coordinates $\{\mc J_{-1}^\alpha\}_\alpha$.  
In particular, both
$\zeta$ and $\mc J$ can be computed recursively by an algebraic algorithm via the
above formula.
\end{prop}

Explicitly, let us denote by $J(\cdot,\cdot)$ and $*$ the flat metric (the first residue pairing) and the product structure on the tangent bundle of $S$, respectively.  For simplicity, let us denote by $t_1,\cdots,t_\mu$ the flat coordinate system on $S$ and by $\partial_{t_1},\cdots,\partial_{t_\mu}$ their partial derivatives. Then, as a consequence of the flat structure, the following  3-tensor
$$
A(\partial_{t_i},\partial_{t_j},\partial_{t_k})
\ := \ J\big(\partial_{t_i} * \partial_{t_j}, \partial_{t_k}\big)  =J\big(\partial_{t_i} , \partial_{t_j}*\partial_{t_k}\big)  \ \ \in \ \Gamma(S,\mathcal{O}_S) 
\qquad 1\le i,j,k \le \mu
$$
is symmetric in the three variables,  and satisfies the following integrability conditions 

$$\partial_{t_l } A(\partial_{t_i},\partial_{t_j},\partial_{t_k})=\partial_{t_i} A(\partial_{t_l},\partial_{t_j},\partial_{t_k})\qquad \text{ for all }  1\le i,j,k,l\le \mu. 
$$

\noindent
Therefore, there exists a function (formal power series in the flat coordinates) $\mc F_{0,f}^{\rm SG}$ on $S$, called the {\it prepotential}, such that
$$\partial_{t_i}\partial_{t_j}\partial_{t_k} \mc F_{0,f}^{\rm SG}=A(\partial_{t_i},\partial_{t_j},\partial_{t_k})=J\big(\partial_{t_i} * \partial_{t_j}, \partial_{t_k}\big)$$
(where the quadratic terms are normalized to be 0).

\noindent We are enabled to compute the {\it prepotential} $\mc F_{0,f}^{\rm SG}$ of the associated formal Frobenius manifold structure in a perturbative way, for an arbitrary weighed homogeneous singularity. On the other hand, it is shown in \cite{LLSaito} that    the formal power series $\zeta$ is in fact the Taylor series expansion of the associated (analytic) primitive form around the origin $\mathbf{0}\in S$. This explains the geometric origin of the induced (formal) Frobenius manifold structure in the above proposition together with the analyticity of its prepotental $\mc F_{0,f}^{\rm SG}$.

Let us restrict our attention to the case of exceptional unimodular singularities now.  Originally,  there are 14 exceptional unimodular singularities   by Arnold \cite{Arnold-strangduality}, which are one parameter families of singularities with three variables. Each family contains a weighted homogenous singularity characterized by the existence of only one negative degree but no zero-degree deformation parameter \cite{Saito-exceptional}. Hence  in the present note,  by  exceptional unimodular singularities, we mean the weighted homogeneous polynomials in these one parameter families, which are given in Table \ref{tab-exceptional-singularities}.

   \begin{center}
    \begin{table}[h]
        \caption{\label{tab-exceptional-singularities}  Exceptional unimodular singularities
    }
        \begin{tabular}{|c|c||c|c||c|c||c|c|}
       \hline
         & Polynomial &  & Polynomial &  & Polynomial &   & Polynomial   \\ \hline \hline

    $E_{12}$ & $x^3+y^7$ &  $W_{12}$ & $x^4+y^5$ &$U_{12}$ & $x^3+y^3+z^4$ & &    \\              \hline
       $Q_{12}$ & $ {x^2y+xy^3+z^3}$  &  $Z_{12}$ & $x^3y+y^4x$ & $S_{12}$ & $x^2y+y^2z+z^3x$   &&  \\              \hline
   $E_{14}$ & $ {x^2+xy^4+z^3}$  &  $E_{13}$ & $x^3+xy^5$  & $Z_{13}$  & $ {x^2+xy^3+yz^3}$&   $W_{13}$ & $ {x^2+xy^2+yz^4}$   \\              \hline
       $Q_{10}$ & $x^2y+y^4+z^3$ & $Z_{11}$ & $x^3y+y^5$    &  $Q_{11}$ & $x^2y+y^3z+z^3$ &$S_{11}$ & $x^2y+y^2z+z^4$     \\              \hline

     \end{tabular}
    \end{table}
     \end{center}

There is a partial classification \cite{Saito-exceptional} of weighted homogeneous polynomial with isolated singularity by using the central charge
  $$ \hat c_f:=\sum_{i=1}^n (1-2q_i).$$
  The case $\hat c_f\leq 1$ is characterized as ADE-singularities if $\hat c_f<1$, or simple elliptic singularities if $\hat c_f=1$.
  The first examples of $\hat c_f>1$ are the    exceptional unimodular singularities, the central charge of which  are   listed in Table \ref{tab-centralcharge} by direct calculations.
     \begin{center}
    \begin{table}[h]
        \caption{\label{tab-centralcharge}
    }

 \begin{tabular}{|c||c|c|c|c|c|c|c|c|c|c|c|c|c|c|}
   \hline
  Type  &$E_{12}$  & $E_{13}$ & $E_{14}$ & $Z_{11}$ & $Z_{12}$ & $Z_{13}$ &$W_{12}$  &$W_{13}$  & $Q_{10}$ & $Q_{11}$ & $Q_{12}$ &$S_{11}$  & $S_{12}$ & $U_{12}$ \\ \hline\hline
  $\hat c_f$  & ${22\over 21}$ & ${16\over 15}$ & ${13\over 12}$ & ${16\over 15}$ & ${12\over 11}$  & ${10\over 9}$ & ${11\over 10}$ & ${9\over 8}$ & ${13\over 12}$ & ${10\over 9}$ & ${17\over 15}$ & ${9\over 8}$ & ${15\over 13}$  & ${7\over 6}$  \\
   \hline
 \end{tabular}
 \end{table}
 \end{center}

For the 14 singularities $f$,  the good basis is   already known to be unique \cite{Hertling-classifyingspace, LLSaito,Saito-uniqueness},
and is simply given by a basis of Jacobian algebra $\Jac(f)$. Hence,
  the primitive form  is unique  (up to a nonzero scalar).
By Proposition \ref{prop-pert}, we can obtain the data on LG B-model at genus zero in a perturbative way, and in particular we can calculate the four-point function $\mc F_0^{(4)}$ (that is, the degree 4 terms of the prepotential $\mc F_{0,f}^{\rm SG}$  with respect to the flat coordinate system) of the Frobenius manifold structure associated to the primitive form.
For instance for $U_{12}$-singularity, $f=x^3+y^3+z^4$, we let $ \{\phi_i\}_i= \{1,z,x,y,z^2,x z,y z,x y,x z^2,y z^2,$  
$x y z,x y z^2\}$.
By direct calculations, we obtain the four-point function in flat coordinates $(t_1, \cdots, t_{12})$ with respect to the primitive form   $\zeta = dxdydz + O(\mathbf{s})$,
\begin{align*} -\mc F^{(4)}_0&=   \frac{1}{8} t_5^2 t_6 t_7+\frac{1}{6} t_3 t_6^2 t_8+\frac{1}{6} t_4 t_7^2 t_8+\frac{1}{4} t_2 t_5 t_7 t_9
                                          +\frac{1}{6} t_3^2 t_8 t_9+\frac{1}{4}t_2 t_5 t_6 t_{10}+\frac{1}{6} t_4^2 t_8 t_{10}\\
                     &\quad +\frac{1}{8} t_2^2 t_9 t_{10}+\frac{1}{8} t_2 t_5^2 t_{11}+\frac{1}{6} t_3^2 t_6 t_{11}
                           +\frac{1}{6}t_4^2 t_7 t_{11}\fbox{$+\frac{1}{18} t_3^3 t_{12}$}\fbox{$+\frac{1}{18} t_4^3 t_{12}$}\fbox{$+\frac{1}{8} t_2^2 t_5 t_{12}$}
\end{align*}
Here we make  boxes for the last three monomials, which will be   compared with the data on the LG $A$-side, studied in the next section.


\section{Mirror construction and the Fan-Jarvis-Ruan-Witten theory}
For the mirror symmetry purpose, we restrict our singularity into an \emph{invertible polynomial}, where the number of variables is the same as the number of monomials in the polynomial. We consider a pair $(W,G)$, where $W$ is an invertible polynomial with $n$ variables $x_1,\cdots,x_n$ and has no monomials of the form $x_ix_j$ for $i\not= j$.
By rescaling the variables, we can always write this polynomial by 
$$W=\sum_{i=1}^n\prod_{j=1}^n x_j^{a_{ij}}.$$
The matrix $E_W:=(a_{ij})_{n\times n}$ of exponents is called the \emph{exponent matrix} of $W$. Let us use ${\rm Aut}(W)$ to denote the \emph{group of diagonal symmetries of} $W$, 
$${\rm Aut}(W):=\left\{{\rm diag}(\lambda_1,\cdots,\lambda_n)\mid W(\lambda_1x_1,\cdots,\lambda_n x_N)=W(x_1,\cdots,x_n), \lambda_i\in\mathbb{C}^*\right\}.$$ 
Then $G$ is a subgroup in ${\rm Aut}(W)$ containing 
$$J_W:={\rm diag}\left(\exp(2\pi\sqrt{-1}q_1),\cdots,\exp(2\pi\sqrt{-1}q_n)\right),$$ 
with $q_1,\cdots,q_n$ are the weights of variables in $W$. 
Berglund and H\"ubsch constructed a mirror polynomial $W^T$ \cite{BH} by taking
$$W^T=\sum_{i=1}^n\prod_{j=1}^n x_j^{a_{ji}}$$
where $E_{W^T}$ is the transpose matrix of $E_W$.
In general, the LG-LG mirror symmetry relates the pair $(W,G)$ to a mirror pair $(W^T,G^T)$, where $G^T$ is constructed by \cite{BHe,K}.

In particular, if $G={\rm Aut}(W)$, then $G^T$ is the group with only an identity element. A LG-LG mirror symmetry conjecture can be formulated as the equivalence of Frobenius manifold structure associated with the primitive form theory of $W^T$ and that associated with the genus-$0$ \emph{Fan-Jarvis-Ruan-Witten theory} (FJRW) theory of $(W, G={\rm Aut}(W))$.

The FJRW theory is introduced by Fan, Jarvis and Ruan in a series of papers \cite{FJR, FJR2}, based on a proposal of Witten \cite{W}. The theory works for the pair $(W,G)$ in general, where $W$ is a weighted homogenous polynomial which has an isolated critical point at the origin and $G$ is a subgroup in ${\rm Aut}(W)$.
The theory also requires that $G$ contains $J_W$.
For technical reasons, $W$ does not contain any monomial term $xy$.  In the present note, we will focus only on the case $G={\rm Aut}(W)$.

For a pair $(W,{\rm Aut}(W))$, there is an FJRW state space $H_{W}$ which collects all ${\rm Aut}(W)$-invariant part of middle dimensional Lefschetz thimble on the fixed locus of each group element $\gamma$ in ${\rm Aut}(W)$,
$$H_{W}:=
\bigoplus_{\gamma\in {\rm Aut}(W)}H^{\rm mid}({\rm Fix}(\gamma);W^{\infty}_{\gamma};\mathbb{C})^{{\rm Aut}(W)}.$$
Here $W^{\infty}_{\gamma}$ is the preimage of $[M,\infty)$, for $M\gg0$, under the real part of $W$ restricted on the fixed locus ${\rm Fix}(\gamma)$.

Fan, Jarvis and Ruan \cite{FJR, FJR2} studied the space of solutions of Witten equations for $W$ $$
\frac{\partial u_i}{\partial \overline{z}} +\overline{\partial_i W}(u_1,\cdots,u_n)=0, \ \ i=1,\cdots,n
$$
where $z$ is a local coordinate of the curve in consideration (but not the formal variable in primitive form theory) 
and $u_i$ ($1\leq i\leq n$) is a section of a line bundle $L_i$ with suitable degrees over the curve (for algebraic construction, see \cite{CLL,PV}),  
and constructed a \emph{cohomological field theory} (in the sense of Kontsevich-Manin \cite{KM}) $\{\Lambda^W_{g,k}:(H_W)^{\otimes k}\rightarrow H^*(\overline{\mathcal{M}}_{g,k},\mathbb{C})\}$ on moduli space of stable curves $\overline{\mathcal{M}}_{g,k}$. As a consequence, this gives the \emph{FJRW} invariants
\begin{equation}\label{FJRW inv}
\LD\alpha_1\psi_1^{\ell_1},\dots,\alpha_k\psi_k^{\ell_k}\RD_{g,k}^W=\int_{\overline{\mathcal{M}}_{g,k}}\Lambda_{g,k}^W(\alpha_1,\dots,\alpha_k)\prod_{j=1}^{k}\psi_j^{\ell_j}, \quad \alpha_j\in H_{W}.
\end{equation}
Here $\psi_j$ is the $j$-th psi-class on $\overline{\mathcal{M}}_{g,k}.$
The genus-$0$ invariants without $\psi$-class involved give a formal Frobenius manifold structure on $H_W$.
The prepotential of this formal Frobenius manifold is
$$\mathcal{F}^{\rm FJRW}_{0,W}=\sum_{\begin{subarray}{l}
        k\geq3
      \end{subarray}}\frac{1}{k!}\LD\mbf{t}_0,\dots,\mbf{t}_0\RD_{0,k}^{W}, \quad
\mbf{t}_0=\sum_{j=1}^{\mu}t_{0, \alpha_j}\,\alpha_j.$$
It is a formal power series of $t_{0, \alpha_j}, j=1,\cdots,\mu$. More generally, the \emph{{\rm FJRW} total ancestor potential} $\mathscr{A}_{W}^{\rm FJRW}$ is defined to be
$$
\mathscr{A}_{W}^{\rm FJRW}
=\exp\left(\sum_{g\geq0}\hbar^{g-1}\sum_{\begin{subarray}{l}
        k\geq0
      \end{subarray}}\frac{1}{k!}\LD\mbf{t}(\psi_1)+\psi_1,\dots,\mbf{t}(\psi_k)+\psi_k\RD_{g,k}^{W}\right).
$$
Here $\mbf{t}(z)=\sum_{m\geq0}\sum_{j=1}^{\mu}t_{m, \alpha_j}\,\alpha_j\,z^m.$

\section{Mirror symmetry for exceptional unimodular singularities}
In \cite{LLSS}, the following isomorphism between two types of Frobenius manifolds is proven.
\begin{thm}\label{g=0-mirror}
 Let $W^T$ be  one of the 14 exceptional unimodular singularities in Table \ref{tab-exceptional-singularities}. There exists a mirror map $\Psi:\Jac(W^T)\cong H_{W}$, which induces an equality
\begin{equation}\label{thm-0}
\mathcal{F}_{0,W^T}^{\rm SG}\ \ = \ \ \mathcal{F}_{0,W}^{\rm FJRW}.
\end{equation}
\end{thm}
The mirror map $\Psi:\Jac(W^T)\to H_{W}$ is constructed by Krawitz \cite{K} and proven that it is a ring isomorphism under a  technical condition that $W$ (in the FJRW side)  is not allowed to be a chain type polynomial with one weight $1/2$. For exceptional unimodular singularities, this condition excludes two examples, $W^T=x^2 y + y^3 z + z^3 (Q_{11}), x^2 y + y^2 z + z^4 (S_{11})$. However, in \cite{LLSS},  the technical condition is removed by using the Getzler's relation in $\overline{\mathcal{M}}_{1,4}$ \cite{Get}.   

The proof of Theorem \ref{g=0-mirror} mainly uses the axioms of cohomological field theories, in particular, the Witten-Dijkgraaf-Verlinde-Verlinde (WDVV) equations. Combined with the special properties on the weights of the exceptional unimodular singularities, it was proved in \cite{LLSS} that both $\mathcal{F}_{0,W^T}^{\rm SG} $ and $\mathcal{F}_{0,W}^{\rm FJRW}$ are determined by the underlying ring isomorphism and a few initial invariants $\LD \cdots \RD_{0,4}$. The invariants on the primitive form theory side can be calculated by the perturbative formula. On the other hand, again by some WDVV equations, the invariants on the FJRW side can be reduced to invariants which can be calculated by the so-called orbifold-Grothendieck-Riemann-Roch formula. Under the mirror map which identifies the deformation parameter space in primitive form side to the FJRW state space together with the ring structure and the inner product, the invariants on both sides are identified up to a scale $-1$. 
Then, by rescaling mirror map appropriately, we obtain the desired equality \eqref{thm-0}.

This equality of the pre-potentials in genus 0 is lifted to  the equality of higher genus potentials as follows.
For the generic point $\bf{s} \in S$ in the universal unfolding $F$ of $W^T$, the $F(\bf{x},\bf{s})$ is a Morse function in $\bf{x}$ so that its Jacobian ring is a direct sum of the one dimensional algebra $\C$. That is, after the Kodaira-Spencer map identification, the Frobenius algebra structure on the tangent space of $S$ at $\bf{s}$ is {\it semi-simple}. Such a generic point is called \emph{semisimple}. Givental defined a \emph{total ancestor potential} (or a higher genus formula) \cite{G2} using only the genus zero data near the generic point together with the knowledge of the Witten-Kontsevich tau-function. The later is just also called the \emph{total ancestor potential} of the Gromov-Witten theory with the target being a point. Teleman \cite{T} proved that this higher genus formula in a cohomological field theory is uniquely determined by the underlying Frobenius manifold at the semisimple point. The origin in the universal unfolding space $S$ is not semisimple, however,  Givental's formula can be uniquely extended to $\mathscr{A}_{W^T}^{\rm SG}$ at the origin by Milanov \cite{M} (see also Coates-Iritani \cite{CI}). The uniqueness of the extension will upgrade Theorem \ref{g=0-mirror} to an identity of higher genus potential function:
$$\mathscr{A}_{W^T}^{\rm SG} \ \ = \ \ \mathscr{A}_{W}^{\rm FJRW}.$$
This completes a proof of LG-LG mirror symmetry.

 \bigskip

\noindent \textbf{Acknowledgement}: 
The  author is partially supported by  JSPS Grant-in-Aid for Scientific Research (A) No. 25247004.

\end{document}